\documentclass[12pt]{amsart}
\usepackage{latexsym}
\usepackage{amssymb}
\usepackage{amsthm}
\usepackage{latexsym}
\usepackage{longtable}
\usepackage{epsfig}
\usepackage{amsmath}
\usepackage{hhline}
\usepackage{color}

\usepackage{listings}

\newtheorem{theorem}{Theorem}

\newtheorem{corollary}{Corollary}

\def\Tr{\mathrm{Tr}}

\usepackage{xcolor}

\title[Partial augmentations of elements]{On partial augmentations  of elements  in  integral group rings}

\author[Bovdi]{Victor Bovdi}
\address{Department of Mathematical Sciences, UAEU, Al-Ain, United Arab Emirates}
\email{vbovdi@gmail.com}
\author[Mar\'oti]{Attila Mar\'oti}
\address{Alfr\'ed R\'enyi Institute of Mathematics,
E\"otv\"os Lor\'and Research Network, P.O. Box 127, H-1364, Budapest, Hungary}
\email{maroti.attila@renyi.hu}

\begin{document}
\maketitle

\begin{abstract}
Inner relations are derived between partial augmentations of  certain elements (units or idempotents)  in   group rings.
\end{abstract}
\section{Introduction}

Let $KG$ be the group ring of a group $G$ over a commutative ring  $K$ with identity. Let $U(KG)$ be  the group of units of  $KG$.
The subgroup
\[
V(KG)=\Big\{\sum_{g \in G} \alpha_{g}g\in U(KG)\quad \mid\quad  \alpha_g\in K,  \quad \sum_{g \in G} \alpha_{g}=1 \Big\}
\]
of  $U(KG)$  is called the {\it normalized group of units} of $KG$. It is easy to see that if $U(K)$ denotes  the  group of units of the ring $K$, then
\[
U(KG) = V(KG) \times  U(K)
\]
and that $G$ is a  subgroup of  $V(KG)$.

For  $g\in G$ let $g^G$ denote the conjugacy  class of $g$ in $G$.  Let  $u=\sum_{g \in G} \alpha_{g} g\in KG$. For $y \in G$ let $\nu_{y}(u)=\sum_{g \in y^{G}} \alpha_{g}$ be the {\it  partial augmentation} of $u$ with respect to $y$. Observe that $\nu_{x}(u)$ is the same for all $x\in y^G$.

The element $\Tr^{(n)}(u)=\sum_{g \in G\{n\}} \alpha_{g} \in K$ is called the $n^{th}$ {\it generalized trace } of the element $u$  (see \cite[p.\,2932]{Artamonov_Bovdi}), where  $G\{n\}$ is the  set of elements   of order $p^n$  of  $G$ where   $n$ is a non-negative integer and $p$ is a prime. Clearly,   $\Tr^{(0)}(u)$ coincides with  $\nu_{1}(u)=\alpha_1$ of $u\in KG$.

Let $K=\mathbb{Z}$, the ring of integers. Let  $u=\sum_{g \in G} \alpha_{g} g\in V(\mathbb{Z}G)$ be  a torsion unit, that is, an element of finite  order $|u|$.  There are  several  connections between  $|u|$,  the partial augmentations $\nu_{g}(u)$ ($g\in  G$) and  $\Tr^{(i)}(u)$ for $i=0,1,\ldots, |u|$. Such a relationship was  first  obtained by Higman and Berman (see \cite[p.\,2932]{Artamonov_Bovdi} or \cite{Sandling}), namely  that $\nu_1(u)=0$ for a finite group $G$. More generally, it is also a consequence of the Higman-Berman Theorem  that    $\nu_{g}(u)=0$ for every central element $g$ of a finite group $G$. The Higman-Berman Theorem  was extended for arbitrary groups $G$ by Bass and Bovdi  (see \cite[Fact 1.2, p.\,2932]{Artamonov_Bovdi} or  \cite[Proposition 8.14, p.\,185]{Bass} and \cite{Bovdi_trace}). Note that it is still an open question whether  $\nu_{g}(u)=0$ for every  central element $g$ of an arbitrary group $G$?

The {\it spectrum}  of a group is the set of orders of its torsion elements. A main unsolved problem in the theory of integral group rings    is the {\it Spectrum Problem $\mathrm{(SP)}$} which says that the spectra of $G$ and $V(\mathbb{Z}G)$ coincide.
A stronger version of $\mathrm{SP}$ was the {\it Zassenhaus Conjecture} $\mathrm{(ZC)}$, which says that for a finite group $G$ each torsion unit of $V(\mathbb{Z}G)$ is rationally conjugate to an element of $G$. The $\mathrm{ZC}$ can also be  reformulated in terms of conditions on $\nu_{g}(u)$  for each torsion unit $u\in V(\mathbb{Z}G)$. A historical overview  of this topic may be found in the survey \cite{Margolis_del_Rio_survey}.

For certain finite groups $G$, the cornerstone for solving the $\mathrm{ZC}$  is  the  so-called  {\it Luthar-Passi method}  introduced in \cite{Luthar_Passi}.
Together with results such as \cite[Proposition 5]{Hales_Luthar_Passi},  \cite[Proposition 3.1]{Hertweck_I},
%\cite{Hertweck_II},
\cite[Proposition 2.2]{Hertweck_III}, \cite{Cohn_Livingstone} and $(p,q)$-character theory  from \cite{Bovdi_Konovalov_Studia}, the Luthar-Passi method provides   ZC for certain groups $G$ (see \cite{Margolis_del_Rio_survey}) as well as    a   counterexample to ZC  (see  \cite{Eisele_Margolis_conterexample}).

After the negative solution of the ZC a question asked  by Bovdi (see \cite[Fact 1.5, p.\,2932]{Artamonov_Bovdi}) is gaining relevance.  Is it true that if $u$ is a torsion unit of $\mathbb{Z}G$ of order $p^n$ where $p$ is a prime and $n$ is a positive integer,  then  $\Tr^{(i)}(u)=0$  for all $i<n$ and $\Tr^{(n)}(u)=1$?

Note that  the above  methods work exclusively only when $G$ is  finite.
With the exception of the Bass-Bovdi Theorem, there is no result which gives a restriction for  $\nu_{g}(u)$  where  $G$ is an infinite group and $u$ is a torsion unit.

Recall that the M\"obius function $\mu$ is defined on the set of positive integers as follows: $\mu(1) = 1$, $\mu(n) = 0$ if $n$ is divisible by the square of a prime, and $\mu(n) = {(-1)}^{\ell}$ if $n = \prod_{i=1}^{\ell} p_{i}$ where $p_{1}, \ldots , p_{\ell}$ are distinct primes.

Our first  result is a new relation between  partial augmentations of a torsion unit of  $\mathbb{Z}G$ where $G$ is  a finite  group.

\begin{theorem}\label{Thm::1}
Let $u\in V(\mathbb{Z}G)$ be a torsion unit of the integral group ring $\mathbb{Z}G$ of a  finite group $G$. Let $k,\; n$ be  positive integers such that $k$ is coprime to  the  exponent of $G$. If $n$ and $k$  are both congruent to $1$ modulo  $|u|$, then for every $s\in G$  we have
\begin{equation}\label{Eq::1}
\begin{split}
\nu_{s}(u) =  \underset{r \mid t \mid n}{\sum} \mu(r) \cdot \Big( \underset{y^{(knr)/t} = x^{k} \sim s}{\underset{x^{G} , \ \exists y \in G:}{\sum}}  \nu_{x}(u)\Big).
\end{split}
\end{equation}
\end{theorem}
\bigskip
\noindent
Formula \eqref{Eq::9}, which is  part of the proof of Theorem \ref{Thm::1},  may be of independent interest.  The proof of Theorem \ref{Thm::1} also depends on the following result   in which  $G$ is not necessarily a finite group and $u$ is not necessarily a unit.

\newpage
\begin{theorem}\label{Thm::2}
Let $u$ be an element  of the integral group ring $\mathbb{Z}G$ of a  group $G$. Let $p$ be a prime and $q = q'\cdot m$ a positive integer such that $m$ is the $p$-part of $q$ and $q'$ is not divisible by $p$. For every $s\in G$ we have
\begin{equation}\label{Eq::2}
\begin{split}
\nu_{s}(u^{q}) \equiv  \underset{r \mid t \mid q'}{\sum} \mu(r) \cdot \Big( \underset{y^{\frac{qr}{t}} = x^{m} \sim s}{\underset{x^{G} , \ \exists y \in G:}{\sum}}  \nu_{x}(u^{q'}) \Big)
\pmod p.
\end{split}
\end{equation}
\end{theorem}

In the special case when  $G$ is a finite group and $u\in \mathbb{Z}G$ is a torsion unit   the main result of Wagner (see \cite{Wagner}) could be compared with  our Theorem \ref{Thm::2}.

Note that Theorem \ref{Thm::2} may be applied to the case when $u$ is a nilpotent element of $\mathbb{Z}G$ with nilpotency index larger than $q'$.

Let $G$ be  a finite group. Let $\mathbb{Q}$ and $\mathbb{C}$ be the fields of rational and complex numbers respectively.   If  $e$ is   an idempotent of $\mathbb{C}G$, then   $\nu_1(e)\in  \mathbb{Q}$ and   $0 < \nu_1(e) < 1$ unless $e\in\{ 0, 1\}$ (see \cite{Zalesskii}). Furthermore,\quad  $|\nu_g(e)|^2 \leq |g^G|\cdot  \nu_1(e)$\quad  (see \cite[Theorem 2, p.\,208]{Weiss_Idempotent}) \; and \quad $\sum_{i=1}^m{|\nu_i(e)|^2}/{|a_i^G|}\leq 1$, \quad where $\{a_1,\ldots, a_m\}$ is a set of representatives of the conjugacy classes of $G$ (see \cite[Corollary 2.6, p.\,2330]{Hales_Luthar_Passi}).

A consequence of Theorem \ref{Thm::2} is    a new relation between the partial augmentations of an  idempotent  in  $\mathbb{Q}G$ where  $G$ is  an arbitrary  group.

\begin{corollary}\label{C::1}
Let $e$ be an idempotent  of  $\mathbb{Q}G$ of a  group $G$.	Let $\beta\in \mathbb{Z}$ such that $u=\beta e\in \mathbb{Z}G$. Let $p$ be a prime and $q = q'\cdot m$ a positive integer such that $m$ is the $p$-part of $q$ and $q'$ is not divisible by $p$. If $p$ does not divide $\beta$, then for every $s\in G$ we have
\begin{equation}\label{Eq::3}
	\begin{split}
\nu_{s}(u) \equiv  		 \underset{r \mid t \mid q'}{\sum} \mu(r) \cdot \Big( \underset{y^{\frac{qr}{t}} = x^{m} \sim s}{\underset{x^{G} , \ \exists y \in G:}{\sum}}  \nu_{x}(u) \Big)  \pmod p.
	\end{split}
\end{equation}
Moreover, if $G$ is finite and \; $p>4q'\cdot |\beta|\cdot |G|^{3/2}$, then in \eqref{Eq::3} equality holds.
\end{corollary}

\section{Proofs}

\begin{proof}[\underline{Proof of Theorem \ref{Thm::2}}]  For elements $x$ and $y$ in $G$ we write $x \sim y$ if $x$ is conjugate to $y$.

Let $s \in G$. We wish to give an expression for $\nu_{s}(u^{q})$. We need some notation.

Consider the set $\mathcal{K} = \{ (g_{1}, \ldots , g_{q}) \in G^{q} \mid g_{1} \cdots g_{q} \sim s \}$. There is a permutation $\pi$ acting on $\mathcal{K}$ by sending $(g_{1}, g_{2}, \ldots , g_{q}) \in \mathcal{K}$ to $(g_{2}, \ldots , g_{q}, g_{1}) \in \mathcal{K}$. Let $t$ be a positive divisor of $q$. Let the union of those $\langle \pi \rangle$-orbits on $\mathcal{K}$ which have lengths dividing $t$ be denoted by
\[
\mathcal{K}_{t} = \{ (g_{1}, \ldots , g_{q}) \in \mathcal{K} \mid g_{i+t} = g_{i} \ \mathrm{for \ every} \ i \ \mathrm{with} \ 1 \leq i \leq q-t \}
\]
and let the union of orbits length $t$  on $\mathcal{K}$ be $\mathcal{K}^{*}_{t}$. Observe that $\mathcal{K} = \mathcal{K}_{q}$ and $\mathcal{K}_{t} = \cup_{r \mid t} \mathcal{K}^{*}_{r}$.

Write $u=\sum_{g\in G}\alpha_gg\in \mathbb{Z}G$. It is easy to see that
\begin{equation}\label{Eq::4}
\nu_{s}(u^{q}) = \sum_{(g_{1}, \ldots , g_{q}) \in \mathcal{K}} \prod_{j=1}^{q} \alpha_{g_{j}} = \sum_{t \mid q} \sum_{(g_{1}, \ldots , g_{q}) \in \mathcal{K}^{*}_{t}} \prod_{j=1}^{q} \alpha_{g_{j}}.
\end{equation}
Since $\mathcal{K}^{*}_{t}$ is the union of all $\langle \pi \rangle$-orbits of length exactly $t$, the multiplicity of each summand in the sum $\sum_{(g_{1}, \ldots , g_{q}) \in \mathcal{K}^{*}_{t}} \prod_{j=1}^{q} \alpha_{g_{j}}$ is divisible by $t$. Thus (\ref{Eq::4}) provides
\begin{equation}\label{Eq::5}
\begin{split}
\nu_{s}(u^{q}) &\equiv \underset{p \nmid t}{\underset{t \mid q}{\sum}} \sum_{(g_{1}, \ldots , g_{q}) \in \mathcal{K}^{*}_{t}} \prod_{j=1}^{q} \alpha_{g_{j}}  \equiv \underset{t \mid q'}{\sum} \sum_{(g_{1}, \ldots , g_{q}) \in \mathcal{K}^{*}_{t}} \prod_{j=1}^{q} \alpha_{g_{j}} \\& \equiv \sum_{t \mid q'} \sum_{(g_{1}, \ldots , g_{q}) \in \mathcal{K}^{*}_{t}} \Big( \prod_{j=1}^{q'} \alpha_{g_{j}} \Big)^{m}
\equiv \sum_{t \mid q'} \sum_{(g_{1}, \ldots , g_{q}) \in \mathcal{K}^{*}_{t}}  \prod_{j=1}^{q'} \alpha_{g_{j}}
\pmod p.
\end{split}
\end{equation}
\noindent
If $f_{1}$ and $f_{2}$ are two functions from $\mathbb{Z}$ to $\mathbb{Z}$ such that\;  $f_{1}(t) = \sum_{r \mid t} f_{2}(r)$, then\quad  $f_{2}(t) = \sum_{r \mid t} \mu(r) f_{1}(t/r)$. This is the M\"obius inversion formula (see \cite[Theorem 2.9,p.\,32]{Apostol}).

For positive integers $t$ and $r$, put
\[
f_{1}(t) = \sum_{(g_{1}, \ldots , g_{q}) \in \mathcal{K}_{t}} \prod_{j=1}^{q'} \alpha_{g_{j}}\quad\text{and}\quad  f_{2}(r) = \sum_{(g_{1}, \ldots , g_{q}) \in \mathcal{K}^{*}_{r}} \prod_{j=1}^{q'} \alpha_{g_{j}}.
\]
The M\"obius inversion formula then yields
\begin{equation}\label{Eq::6}
\sum_{(g_{1}, \ldots , g_{q}) \in \mathcal{K}^{*}_{t}} \prod_{j=1}^{q'} \alpha_{g_{j}} = \sum_{r \mid t} \mu(r) \sum_{(g_{1}, \ldots , g_{q}) \in \mathcal{K}_{t/r}} \prod_{j=1}^{q'} \alpha_{g_{j}}.
\end{equation}
Formulas (\ref{Eq::5}) and (\ref{Eq::6}) yield
\[
\begin{split}
\quad &\nu_{s}(u^{q}) \equiv \\
&\equiv \underset{t \mid q'}{\sum} \Big( \sum_{r \mid t} \mu(r) \sum_{(g_{1}, \ldots , g_{q}) \in \mathcal{K}_{t/r}} \prod_{j=1}^{q'} \alpha_{g_{j}} \Big) \equiv \underset{r \mid t \mid q'}{\sum} \mu(r) \cdot \Big( \sum_{(g_{1}, \ldots , g_{q}) \in \mathcal{K}_{t/r}} \prod_{j=1}^{q'} \alpha_{g_{j}} \Big) \\
&\equiv \underset{r \mid t \mid q'}{\sum} \mu(r) \cdot \Big( \underset{(g_{1} \cdots g_{t/r})^{{qr}/{t}} = (g_{1} \cdots g_{q'})^{m} \sim s}{\underset{(g_{1}, \ldots , g_{q'}) \in G^{q'}}{\sum}} \prod_{j=1}^{q'} \alpha_{g_{j}} \Big)\\
&\equiv \underset{r \mid t \mid q'}{\sum} \mu(r) \cdot \Big( \underset{y^{{qr}/{t}} = x^{m} \sim s}{\underset{x^{G} , \ \exists y \in G:}{\sum}} \nu_{x}(u^{q'}) \Big)
\pmod p.
\end{split}
\]
\end{proof}

\begin{proof}[\underline{Proof of Theorem \ref{Thm::1}}]
Let $s,p,q,q'$ and $m$ be as in Theorem  \ref{Thm::2}.
Let $n = q'$ and $m = p$. By \eqref{Eq::2} of Theorem \ref{Thm::2} we have
\begin{equation}\label{Eq::7}
\begin{split}
\nu_{s}(u^{np}) \equiv  \underset{r \mid t \mid n}{\sum} \mu(r) \cdot \Big( \underset{y^{{npr}/{t}} = x^{p} \sim s}{\underset{x^{G} , \ \exists y \in G:}{\sum}}  \nu_{x}(u^{n}) \Big)
\pmod p.
\end{split}
\end{equation}
Let $k$ be a positive integer  coprime to the exponent $e$ of $G$.
Choose $p$ such that $p \equiv k \pmod e$. There are infinitely many such primes  by Dirichlet's theorem on arithmetic progressions \cite[Chapter 7]{Apostol}.

Since $p \equiv k \pmod e$, in \eqref{Eq::7} we have $y^{\frac{npr}{t}} = y^{\frac{nkr}{t}}$  and $x^{p}=x^{k}$. Moreover, $u^{np}=u^{nk}$ by the Cohn-Livingstone Theorem \cite[Corollary 4.1]{Cohn_Livingstone}. This yields
\begin{equation}\label{Eq::8}
\begin{split}
\nu_{s}(u^{nk}) \equiv  \underset{r \mid t \mid n}{\sum} \mu(r) \cdot \Big( \underset{y^{{nkr}/{t}} = x^{k} \sim s}{\underset{x^{G} , \ \exists y \in G:}{\sum}}  \nu_{x}(u^{n}) \Big)
\pmod p.
\end{split}
\end{equation}
The absolute value of every partial augmentation of $G$ is at
most $\sqrt{|G|}$ (really $\nu_y(x)^2\leq |y^G|$) by
\cite[Corollary 2.3, p.\,2329]{Hales_Luthar_Passi} or \cite{Bovdi_Hertweck}.
The number of summands on the right-hand side of  \eqref{Eq::8} is at most $(2\sqrt{n})^2\cdot|G|$. Choose $p$ such that $p>  (2\sqrt{n})^2\cdot|G|^{3/2}$.  Since both sides of the congruence \eqref{Eq::8} have  absolute value less than $p$,
\begin{equation}\label{Eq::9}
\begin{split}
\nu_{s}(u^{nk}) =  \underset{r \mid t \mid n}{\sum} \mu(r) \cdot \Big( \underset{y^{(knr)/t} = x^{k} \sim s}{\underset{x^{G} , \ \exists y \in G:}{\sum}} \nu_{x}(u^{n})\Big).
\end{split}
\end{equation}
If  $k$ and $n$ are both congruent to $1$ modulo  $|u|$, then we get \eqref{Eq::1}.
\end{proof}

\begin{proof}[\underline{Proof of Corollary \ref{C::1}}] Let $s,p,q,q'$ and $m$ be as in Theorem \ref{Thm::2}. Since $u^{r}=\beta^{r-1}u$, we get  $\nu_{s}(u^{r})=\nu_{s}(\beta^{r-1}u)=\beta^{r-1}\nu_{s}(u)$, where $r\in\{q,q'\}$. Theorem \ref{Thm::2} gives
\[
\begin{split}
\beta^{q-q'}\nu_{s}(u) \equiv  \underset{r \mid t \mid q'}{\sum} \mu(r) \cdot \Big( \underset{y^{{qr}/{t}} = x^{m} \sim s}{\underset{x^{G} , \ \exists y \in G:}{\sum}}  \nu_{x}(u) \Big)
\pmod p.
\end{split}
\]
Congruence \eqref{Eq::3} follows by observing that $\beta^{q-q'}=(\beta^{m})^{q'}\beta^{-q'}\equiv 1\pmod{p}$ since  $m$ is a $p$-power.

The absolute value of the left-hand side of \eqref{Eq::3} is at most $|\beta|\cdot \sqrt{|G|}$ and the absolute value of the right-hand side of \eqref{Eq::3} is at most\;  $(2\sqrt{q'})^2\cdot |G|\cdot |\beta|\cdot \sqrt{|G|}$,\; by \cite[Theorem 2, p.\,208]{Weiss_Idempotent}. If  \; $p>4q'\cdot |\beta|\cdot |G|^{3/2}$, then equality in \eqref{Eq::3} holds.
\end{proof}
%{\bf Statement of competing interests.} We declare that there are no competing interests between authors.

\end{document}